\documentclass[11pt,one side]{amsart}
\usepackage{graphicx}
\usepackage{amsmath}
\usepackage{amsthm}
\usepackage{amsfonts}
\usepackage{amssymb}

\theoremstyle{plain}
\newtheorem{theorem}{Theorem}
\newtheorem{proposition}{Proposition}
\newtheorem{corollary}{Corollary}
\theoremstyle{definition}
\newtheorem{remark}{Remark}
\theoremstyle{remark}
\newtheorem*{example}{Example}

\newcommand{\al}{\alpha}
\renewcommand{\a}{\alpha}

\newcommand{\Z}{\mathbb Z}

\renewcommand{\u}{\nearrow}
\renewcommand{\d}{\searrow}

\begin{document}

\begin{center}
{\Large L'Hospital-type rules for monotonicity and limits:\\
Discrete case}

\bigskip

{\large Iosif Pinelis}

\medskip

{\em Department of Mathematical Sciences, Michigan Technological University, Houghton, Michigan 49931 USA}

\end{center}


\bigskip

\section{Introduction}
\ 
Let $-\infty\leq a<b\leq\infty$, let $f$ and $g$ be continuously differentiable functions
defined on the interval $(a,b)$, and let $r=f/g$ and $\rho={f'}/{g'}$. 
In {\cite{3}}, general ``rules" for monotonicity patterns, resembling the usual l'Hospital rules for limits, were given.
For example, according to Proposition~1.9 in
{\cite{3}}, one has the following:
if $\rho$ is increasing and $gg'>0$ on $(a,b)$, then $r\d\u$, which means that 
there is some $c$ in $[a,b]$ such that $r$ is decreasing on $(a,c)$ and 
increasing on $(c,b)$. 
In particular, if $c$ is either $a$ or $b$, the result is that $r$ is either increasing or decreasing on the entire interval $(a,b)$. 
If one also knows whether $r$ is increasing or decreasing in a right neighborhood of $a$ and in a left neighborhood of $b$, then one can discriminate with certainty between these three patterns. 
Using such rules, 
one can generally determine ({\cite{3,monthly}}) the monotonicity pattern of $r$ given that of $\rho$, however complicated the latter might be. 

Clearly, these l'Hospital-type rules for monotonicity patterns are helpful wherever the l'Hospital rules for limits are. Moreover, the monotonicity rules apply even outside such contexts, because they do not require that both $f$ and $g$ (or either of them) tend to 0 or $\infty$ at any point. 
In the special case when both $f$ and $g$ vanish at an endpoint of the interval $(a,b)$, l'Hospital-type rules for monotonicity can be found, in different forms and with different proofs, in {\cite{aqvv}--\cite{avv01}}, {\cite{chavel}--\cite{grom2}}, {\cite{pin91}}, and {\cite{2}--\cite{5}}.
In view of what has been said here, it should not be surprising that a very wide variety of applications of these l'Hospital-type rules for monotonicity patterns were given in those papers; see also {\cite{monthly}}.
In the present note, discrete analogues of l'Hospital-type rules both for monotonicity and limits are given. 

\section{Rules}
\newcommand{\be}[2]{\rule{0pt}{12pt}\overline{#1,#2}}
\newcommand{\de}{\Delta}
\renewcommand{\le}{\leqslant}
\renewcommand{\le}{\leqslant}
\renewcommand{\ge}{\geqslant}
 
Let $f:=(f_n\colon n\in\be ab\,)$ and $g:=(g_n\colon n\in\be ab\,)$ be two real sequences or, equivalently, functions defined on an interval $\be ab:=\{n\in\Z\colon a\le n\le b\}$ of integers, where $a$ and $b$ are in $\Z\cup\{-\infty,\infty\}$ 
(so that $\be ab=\emptyset$ if $a>b$, and no infinite endpoint belongs to $\be ab$). Let
$$r:=\frac fg\quad\text{and}\quad\rho:=\frac{\de f}{\de g},$$
where $(\de f)_n:=\de f_n:=f_n-f_{n-1}$ for $n\in\be{a+1}b$, so that the function 
$\de f$ is defined on $\be{a+1}b$ (with $\pm\infty+1:=\pm\infty$). 
It is assumed throughout that $g$ and $\de g$ do not take on the zero value and do not change their respective signs.

\begin{theorem} \label{th}
Suppose that $\rho$ is either nondecreasing or nonincreasing. 
Then the dependence of the monotonicity pattern of $r$ on that of 
$\rho$ (and also on the sign of $g\,\de g$) is given by the following table:\\


\noindent
\parbox[l]{1.3in}
{
Table 1.\hfill
\vspace*{6pt}

\begin{tabular}{|c|c||c|}
\hline
$\rho$	& $g\,\de g$	& $r$ \\ \hline\hline
$\u$		& $>0$	& $\d\u$ \\ \hline
$\d$		& $>0$	& $\u\d$ \\ \hline
$\u$		& $<0$	& $\u\d$ \\ \hline
$\d$		& $<0$	& $\d\u$ \\ \hline 
\end{tabular}
}
\parbox[l]{4.0in}
{
Here, for instance, the statement $r\d\u$ can be taken to mean that there is some $k$ in $\be ab\,\cup\{a,b\}$ such that $r$ is nonincreasing (\,$\d$\,) on $\be ak$ and 
nondecreasing ($\u$) on $\be kb$. 
In particular, if $k=a$ then $r\d\u$ will imply $r\u$ on the entire interval $\be ab$; similarly, if $k=b$ then $r\d\u$ will imply $r\d$ on $\be ab$. 
}

\vspace*{4pt}

\end{theorem}  

\begin{remark} \label{rem:discrim}
To discriminate between these three possibilities ($k=a$, $k=b$, and $a<k<b$) in the case when (say) $a$ and $b$ are finite, it suffices to know whether $r_{a+1}\ge r_a$ and 
$r_{b-1}\ge r_b$; if (say) $b=\infty$, then one may instead want to know the monotocity pattern of $r$ in a neighborhood of $\infty$.
\hfill $\diamondsuit$
\end{remark}   

\begin{proof}[Proof of Theorem~\ref{th}]
Without loss of generality, $a$ and $b$ are finite.
In view of the ``horizontal" and ``vertical" reflections 
$\Z\ni n\leftrightarrow(-n)$ and 
$f\leftrightarrow(-f)$, it suffices to consider only the first line of Table~1, with 
$\rho\u$ and $g\,\de g>0$. 
Next, it suffices to show that there is some $k\in\be ab$ such that $\de r\le0$ on 
$\be{a+1}k$ and $\de r>0$ on 
$\be{k+1}b$. 

Let now $k:=\sup\{n\in\be{a+1}b\colon \de r_n\le0\}$ (here, $\sup\emptyset:=a$). 
Then $\de r>0$ on $\be{k+1}b$. 
If $\de r\le0$ on $\be{a+1}k$, the proof is completed. 
To obtain a contradiction, assume the contrary, that 
$\de r_n>0$ for some $n\in\be{a+1}k$ (so that $k\ne a$). 
Hence, there exists $m:=\max\{n\in\be{a+1}k\colon\de r_n>0\}$.
Then in fact $m\in\be{a+1}{k-1}$; this follows because  
$k\ne a$ and hence, by the definition of $k$, one has $\de r_k\le0$, while 
$\de r_m>0$. 
Now it also follows from the definition of $m$ that $\de r_{m+1}\le0$. 

The key observation is that for all $n\in\be{a+1}b$
\begin{equation} \label{eq:identities}
g_n g_{n-1}\,\de r_n=(\rho_n-r_n)\,g_n\,\de g_n
=(\rho_n-r_{n-1})\,g_{n-1}\,\de g_n.
\end{equation}
Using these identities (with $n=m$ and $n=m+1$) together with the obtained above inequalities 
$\de r_m>0$ and $\de r_{m+1}\le0$, one concludes that
$\rho_m>r_m\ge\rho_{m+1}$, which contradicts the condition that $\rho$ is nondecreasing. 
\end{proof}

\begin{remark} \label{rem:strict}
{\bf (i)}\ 
For the case given by the first line of Table~1, when $\rho\u$ and 
$g\,\de g>0$, 
the above proof shows that there is some $k\in\be ab$ such that $\de r\le0$ on 
$\be{a+1}k$ and $\de r>0$ on 
$\be{k+1}b$. 
Using the horizontal reflection $\Z\ni n\leftrightarrow(-n)$, one can then see that there also exists some $\ell\in\be ab$ such that $\de r<0$ on 
$\be{a+1}\ell$ and $\de r\ge0$ on 
$\be{\ell+1}b$.
Hence, the conclusion $r\d\u$ for this first-line case in Table~1 can actually be understood in a slightly stronger sense: that there are some $k$ and $\ell$ in $\be ab\,\cup\{a,b\}$ such that $r$ is 
(strictly) decreasing on 
$\be a\ell$, constant on $\be\ell k$, and 
increasing on $\be kb$.
Moreover, identities \eqref{eq:identities} show that $r$ equals a constant $C$ on an interval 
$\be\ell k$ of integers only if $\rho=C$ on $\be{\ell+1}k$ (indeed, if $r=C$ on $\be\ell k$, then $\de r=0$ on $\be{\ell+1}k$, and so, 
$\rho_n=r_n=C$ for all $n\in\be{\ell+1}k$, because of the requirement that neither $g$ nor $\de g$ vanish at any point). 

{\bf (ii)}\ It follows that, if $\rho$ is (strictly) increasing or decreasing, then the statement 
$r\d\u$ can be taken to mean that there is some $k$ in 
$\be ab\,\cup\{a,b\}$ such that $r$ is decreasing on $\be ak$ and 
increasing on $\be{k+1}b$. 
For further details on this point, see the example at the end of this note.

{\bf (iii)}\ Similar comments are valid for the other three cases given by the Table~1.
\ \flushright $\diamondsuit$
\end{remark}

\medskip

Proposition~1.9 in {\cite{3}}, mentioned earlier, was in fact a corollary of a general result, Proposition~1.2 in {\cite{3}}, which also contains the special case when both $f$ and $g$ vanish at an endpoint of the interval $(a,b)$ (as pointed out in Remark~1.5 in {\cite{3}}. Here, we shall treat the discrete analogue of that special case separately, as follows. 

\begin{proposition} \label{prop:special}  
Suppose that either (i) $b=\infty$, $f_\infty:=\lim_{n\to\infty}f_n=0$, and $g_\infty=0$ or (ii) $a=-\infty$, $f_{-\infty}:=\lim_{n\to-\infty}f_n=0$, and $g_{-\infty}=0$.  
Suppose also that $\rho$ is nondecreasing or nonincreasing or increasing or decreasing; then $r$ is so, respectively.  
\end{proposition}

\begin{proof}
Consider the case when $b=\infty$, $f_\infty=0$, and $g_\infty=0$. Then the result follows immediately from the identity
$$g_n g_{n-1}\,\de r_n=
\sum_{j=n}^\infty \de g_n\,\de g_j\,(\rho_j-\rho_n)$$
for all $n\in\be{a+1}\infty$.
The other case can now be proved by the horizontal reflection; compare with the proof of Theorem~\ref{th}. 
\end{proof}

``Discrete" analogues of l'Hospital's rules for limits are easily stated and proved, and apparently well known. Let us present them here for the sake of completeness.  

\begin{proposition} \label{prop:limits}  
Suppose that $b=\infty$, and either 
(i) $|g_\infty|=\infty$ or (ii) $f_\infty=g_\infty=0$. Then $r_\infty=\rho_\infty$ provided that the latter limit, $\rho_\infty$, exists.  
\end{proposition}

\begin{proof}
Let us first consider part (i), assuming that $|g_\infty|=\infty$. 
For $m<n$, one has 
$r_n=\dfrac{f_m}{g_n}+r_{m,n}\Big(1-\dfrac{g_m}{g_n}\Big)$, where 
$r_{m,n}:=\dfrac{f_n-f_m}{g_n-g_m}=
\sum\limits_{j=m+1}^n\rho_j\,|\de g_j|\,\Big/\sum\limits_{j=m+1}^n|\de g_j|$, 
since $\de g$ is assumed not to change sign.
It follows that 
$r_{m,n}$ is between the minimum and maximum values of $\rho$ over the interval $\be{m+1}n$. To complete the proof of part (i), it remains to let $n\to\infty$ and then $m\to\infty$. 

Part (ii) is even simpler to prove, since here $r_m=r_{m,\infty}$ for all 
$m$. 
\end{proof}

Note that condition $|f_\infty|=\infty$ is not required in part (i) of Proposition~\ref{prop:limits}. Moreover, the above proof can be obviously adapted (using, say, the Mean-Value Theorem for the ratio corresponding to $r_{m,n}$) to 
l'Hospital's original ``differentiable" case, even though the condition 
$|f_\infty|=\infty$ is traditionally included into the formulation of the corresponding l'Hospital rule for limits.  

\section{Illustrations}
Let $p:=(p_j\colon j\in\Z)$ be a positive sequence. 
Then $p$ is called log-convex on $\be ab$ if 
$q:=\ln p$ is convex on $\be ab$, in the sense that $\de q$ is nondecreasing on 
$\be{a+1}b$, which is equivalent to the condition that the ratio $p_n/p_{n+1}$ be nonincreasing in $n\in\be a{b-1}$. 
The log-concavity of a sequence and the strict versions of log-convexity and log-concavity are defined similarly.  

\begin{corollary} \label{prop:log-con-tail} 
Suppose that $p$ is log-convex or log-concave on $\be a\infty$, and 
$f_n:=\sum_{j=n}^\infty p_j<\infty$ for all $n\in\Z$.   
Then $f$ is, respectively, log-convex or log-concave on $\be a\infty$.
More generally, the same conclusion holds for any natural $k$ if $f=R^k p$, where $R^k p$ is given by the formula
$$(R^k p)_n:=\sum_{j=n}^\infty\binom{j-n+k-1}{j-n} p_j
\quad\text{for all $n\in\Z$.}$$ 
\end{corollary}

\begin{proof}
Let $g_n:=f_{n+1}$ for all $n\in\Z$. Then the first part of Corollary \ref{prop:log-con-tail} follows immediately from Proposition~\ref{prop:special}.
In turn, this yields the second part of the corollary, because it is easy to see that $(R^k\colon k\in\be1\infty\,)$ is a semigroup of operators with 
$(R^1 p)_n=\sum_{j=n}^\infty p_j$ for all $n$; cf. Remark 5 in \cite{pin99}.
\end{proof}

Corollary~\ref{prop:log-con-tail} is essentially well-known. For the ``log-concave" part, see, for example, {\cite{pin99}} (where $k$ was allowed to be any positive real number)
and, for the continuous counterpart in the case $k=1$, 
{\cite{barlow}} and {\cite{perlman}}. 
The ``log-convex" part can also be obtained from the well-known fact that any linear combination with positive coefficients of log-convex functions is log-convex, having also in mind that the log-convexity is preserved under the shift $n\mapsto n+1$. 
	
\begin{corollary} \label{prop:log-con} 
Suppose that $p$ is log-concave on $\be0\infty$, and 
$f_n:=\sum_{j=0}^n p_j$ for all $n\in\be0\infty$.   
Then $f$ is {\em strictly} log-concave on $\be0\infty$.
More generally, the same conclusion holds for any natural $k$ if $f=L^k p$, where $L^k p$ is given by the formula
$$(L^k p)_n:=\sum_{j=0}^n\binom{n-j+k-1}{n-j} p_j
\quad\text{for all $n\in\be0\infty$.}$$
\end{corollary}

\begin{proof}
Let again $g_n:=f_{n+1}$ for all $n\in\be0\infty$. Then the first part of Corollary~\ref{prop:log-con} follows immediately from part (i) of Remark~\ref{rem:strict} (since
$g_0 g_1\,\de r_1>p_1^2-p_0\,p_2\ge0$, and so, $r_1>r_0$).
In turn, this yields the second part of the corollary, since $(L^k\colon k\in\be1\infty\,)$ is a semigroup with $(L^1 p)_n=\sum_{j=0}^n p_j$ for all $n$. (One can note that $(L^k p)_n=(T^{-1}R^k Tp)_n$ for all $n\in\be0\infty$, all natural $k$, and all $p$ such that $p=0$ on $\be{-\infty}{-1}$, where $(Tp)_n:=p_{-n}$ for all $n$.)
\end{proof}

However, the ``log-convex" analogue of Corollary~\ref{prop:log-con} does not hold. 
Indeed, if a sequence $p$ is both log-convex and log-concave (that is, geometric) then, by Corollary~\ref{prop:log-con}, $f$ is strictly log-concave and hence not log-convex. 

Let us conclude this paper with another illustration of presented results and methods; note the use of 
Remarks \ref{rem:discrim} and \ref{rem:strict} and identities 
\eqref{eq:identities}. 

\begin{example}
Let sequences $f:=f^{(\al)}$ and $g$ on $\be0\infty$ be given by the formulas
$f_n^{(\al)}:=\al+\sum_{j=0}^n p_j$ and 
$g_n:=\sum_{j=0}^n q_j$
for any $\al\ge0$ and all $n\in\be0\infty$, where $p$ and $q$ are positive sequences such that $\rho=p/q$ is increasing on $\be1\infty$ and $r_0^{(0)}<r_1^{(0)}$, where in turn $r^{(\al)}:=f^{(\al)}/g$. 
(For instance, one can take $p_j=j!$ and $q_j=(j/e)^j$ for all 
$j\in\be0\infty$, assuming that 
$0^0=1$.)
 
Then, by Theorem~\ref{th} (or, rather, by part (ii) of Remark~\ref{rem:strict}; cf. Remark~\ref{rem:discrim}), the sequence 
$r^{(0)}$ is increasing on the entire interval $\be0\infty$. 
Moreover, for every $\al>0$ there is some $k_\al$ in $\be0\infty\,\cup\{\infty\}$ such that $r^{(\al)}$ is decreasing on $\be 0{k_\al}$ and 
increasing on $\be{k_\al+1}\infty$. 

The observed condition that $r^{(0)}$ is increasing on $\be0\infty$
implies $\de r^{(0)}>0$ and hence, by \eqref{eq:identities}, 
$\rho>r^{(0)}$ on $\be1\infty$. On the other hand, 
$r^{(\al)}=r^{(0)}+\al/g$. 
Therefore, for each $k\in\be1\infty$, the number $\al_k:=(\rho_k-r_k^{(0)})\,g_k$ is positive and satisfies the equation $\rho_k=r_k^{(\al_k)}$, so that, by \eqref{eq:identities}, $(\de r^{(\al_k)})_k=0$; that is, $r_k^{(\al_k)}=r_{k-1}^{(\al_k)}$. 
Using again part (ii) of Remark~\ref{rem:strict}, one sees that $r^{(\al_k)}$ is decreasing on $\be 0{k-1}$, constant on $\{k-1,k\}$, and 
increasing on $\be k\infty$.
It follows, in particular, that one cannot replace $\be{k+1}b$ in part (ii) of Remark~\ref{rem:strict} by $\be kb$ (keeping the rest unchanged). 

In view of \eqref{eq:identities}, observe also that $r_k^{(\al)}\ge r_{k-1}^{(\al)} \iff \rho_k\ge r_k^{(\al)} \iff \al\le\al_k$, for all $\al\ge0$, and these equivalences hold if all three non-strict inequalities here replaced by the corresponding strict ones. In particular, one has $\al\le\al_k\implies
r_k^{(\al)}\ge r_{k-1}^{(\al)}\implies r_{k+1}^{(\al)}>r_k^{(\al)}\implies
\al<\al_{k+1}$; the second implication here follows again by part (ii) of Remark~\ref{rem:strict}.
Hence, $\a_k$ is increasing in $k$. 
Moreover, if $\al\in(\al_k,\al_{k+1})$ then $r_{k-1}^{(\al)}>r_k^{(\al)}<r_{k+1}^{(\al)}$, so that $r_n^{(\al)}$ is decreasing in $n\in\be1k$ and increasing in $n\in\be k\infty$; recall that here $k$ was taken to be any natural number. 
 
\vspace*{6pt}
\noindent
\hspace*{-15pt}
\parbox{1.9in}{\includegraphics{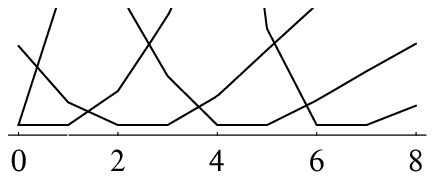} }
\parbox{3.55in}
{
To illustrate, here are parts of the graphs of the linear interpolations of the sequences $n\mapsto r_n^{(\al_k)}/r_k^{(\al_k)}-0.97$ for $k=0,1,3,5,7$  
where $\al_0:=0$ and, as above, \\
$p_j=j!$ and $q_j=(j/e)^j$ with $0^0:=1$.
}
 
\vspace*{6pt}
To visualize the idea of this example, one can imagine a tank with the solution of a liquid in water. Initially, at time $n=0$, the amounts of the liquid and water are $f_0^{(\al)}=\al+p_0$ and 
$g_0=q_0$, respectively, so that the initial relative concentration of the liquid (with respect to water) is $r_0^{(\al)}=(\al+p_0)/q_0$. Suppose that, at each of the time moments $n=1,2,\dots$, the liquid and water are added to the tank in the amounts of $\de f_n=p_n$ and $\de g_n=q_n$, respectively, so that the relative concentration of the liquid in the $n$th addition is $\rho_n=p_n/q_n$ and that in the tank at time $n$ is $r_n^{(\al)}$. 
If $\al$ is large enough, then initially the relative concentration $\rho$ of the liquid in what is added is less than the relative concentration $r^{(\al)}$ of the liquid in the tank, so that $r^{(\al)}$ will be decreasing in time. However, at least in the case when $\rho$ is increasing to $\infty$, $\rho$ will eventually overtake $r^{(\al)}$, and the latter will then be forever increasing. 
\hfill $\diamondsuit$
\end{example}

\bigskip

\end{document}